\documentclass{amsproc}
\usepackage{amssymb}
\usepackage{amsfonts}

\theoremstyle{plain}
\newtheorem{acknowledgement}{Acknowledgement}

\newtheorem{corollary}{Corollary}

\newtheorem{definition}{Definition}

\newtheorem{problem}{Problem}
\newtheorem{proposition}{Proposition}

\newtheorem{theorem}{Theorem}
\numberwithin{equation}{section}
\input{tcilatex}

\begin{document}
\title[Parabolic Stein manifolds]{Parabolic Stein Manifolds}
\author{A. Aytuna}
\address{FENS, SABANCI UNIVERSITY, 34956 TUZLA ISTANBUL,TURKEY}
\email{aytuna@sabanciuniv.edu}
\author{A. Sadullaev}
\address{MATHEMATICS DEPARTMENT , NATIONAL UNIVERSITY OF UZBEKISTAN VUZ
GORODOK, 700174, TASKENT, UZBEKISTAN}
\email{sadullaev@mail.ru}
\date{April 27, 2010}
\thanks{The first author is partially supported by a grant from Sabanci
University\\
The second author is partially supported by Khorezm Mamun Academy, Grant $%
\Phi A\Phi 1\Phi 024.$}

\begin{abstract}
An open Riemann surface is called parabolic in case every bounded
subharmonic function on it reduces to \ a constant. Several authors
introduced seemingly different analogs of this notion for Stein manifolds of
arbitrary dimension. In the first part of this note we compile these notions
of parabolicity and give some immediate relations among these different
definitions. In section 3 we relate some of these notions to the linear
topological type of the Fr\'{e}chet space of analytic functions on the given
manifold. In sections 4\ and 5 we look at some examples and show, \ for
example, that the complement of the zero set of a Weierstrass polynomial
possesses a continuous plurisubharmonic exhaustion function that is maximal
off a compact subset.
\end{abstract}

\maketitle

\section{\textbf{Introduction}}

In the theory of Riemann surfaces, simply connected manifolds, which are
equal to the complex plane are usually called parabolic and\ the ones which
equal to the unit disk are called hyperbolic. Several authors introduced
analogs of these notions for general complex manifolds of arbitrary
dimension in different ways; in terms of triviality (parabolic type) and
non-triviality (hyperbolic type) of the Kobayashi or Caratheodory metrics,
in terms of plurisubharmonic $(psh)$ functions, etc. In some of these
considerations existence of rich family of bounded holomorphic functions
play a significant role.

On the other hand, attempts to generalize Nevanlinna's value distribution
theory to several variables by Stoll, Griffiths, King, et al. produced
notions of "parabolicity" in several complex variables defined by requiring
the existence of certain special plurisubharmonic functions. The common
features of these special plurisubharmonic functions were that they were
exhaustive and maximal outside a compact set.

Following Stoll, we will call an n-dimensional complex manifold $X,$ $%
S-parabolic$\ in case there is a\ plurisubharmonic function $\rho $ on $X$
with the properties:

a) $\left\{ {z\in X:\rho \left( z\right) \leqslant C}\right\} \,\,\subset
\subset X\,\,,\ \forall \,C\in {R}^{+}\,\,$ (i.e. ${\rho }$ is exhaustive),

b) The Monge - Amp\`{e}re operator $\left( {dd^{c}\rho }\right) ^{n}$ is
zero off a compact $\mathrm{K}\subset \subset X$. That is $\rho $ a maximal
plurisubharmonic function outside K.

We note that without the maximality condition b), an exhaustion function $%
\sigma \left( z\right) \in psh\left( X\right) \cap C^{\infty }\left(
X\right) $ always exist for any Stein manifold $X$, because such manifolds
can be properly embedded in $\mathbb{C}_{w}^{2n+1}$ and one can take for $%
\sigma $ the restriction of $\ln \left\vert w\right\vert $ to $X$.

On $S$-parabolic manifolds any bounded above plurisubharmonic function is
constant. In particular, there are no non-constant bounded holomorphic
functions on such manifolds. Complex manifolds, on which every bounded above
plurisubharmonic function reduces to a constant, a characteristic shared by
parabolic open Riemann surfaces and affine-algebraic manifolds, play an
important role in the structure theory of Fr\'{e}chet spaces of analytic
functions on Stein manifolds, especially in finding continuous extension
operators for analytic functions from complex submanifolds (see, \cite{V12},%
\cite{AY2},\cite{AY4}). Such spaces will be called $"parabolic"$ in this
paper.

Special exhaustion functions with certain regularity properties play a key
role in the Nevanlinna's value distribution theory of holomorphic maps $%
f:X\rightarrow \mathrm{P}^{m},$ where $\mathrm{P}^{m}\,-\,m$ dimensional
projective manifold (see.\cite{SA2},\cite{SH},\cite{ST1},\cite{ST2}).\ On
the other hand for manifolds which have a special exhaustion function one
can define extremal Green functions as in the classical case and apply
pluripotential theory techniques to obtain analogs of some classical results
which were proved earlier for $\mathbb{C}^{n}$ (\cite{D}, \cite{ZE1}). In
the special case of an affine algebraic manifolds such a program was
successfully carried out in \cite{ZE2}

The aim of this paper is to state and analyze different notions of
parabolicity, give examples and relate the parabolicity of a Stein manifold $%
X$ with the linear topological properties of the Fr\'{e}chet space of global
analytic functions on $X$.

The organization of the paper is as follows: In section 2, we state and
compare different definitions of parabolicity.\textit{\ }We also bring to
attention, a problem in complex pluripotential theory that arise in this
context. In the third section we relate the notion of parabolicity of a
Stein manifold $X$ with the linear topological type of the Fr\'{e}chet space
,$O(X)$, of analytic functions on $X$. We introduce the notion of \textit{%
tame isomorphism} of $O(X)$ to the space of entire functions (Definition 3)
and show (Theorem 2) that a Stein manifold of dimension $n$ is $S^{\ast }$%
-parabolic if and only if $O(X)$ is tamely isomorphic to the space of entire
functions in $n$ variables. The final section is devoted to some classes of
parabolic manifolds. First we look at complements of the zero sets of entire
functions and show that the complement in $\mathbb{C}^{n}$, of the zero set
of a global Weierstrass polynomial ( algebroidal function), is $S^{\ast }$%
-parabolic. Then we generalize a condition of Demailly for parabolicity and
use it to show that Sibony-Wong manifolds (see section 4 for the definition)
are parabolic.

Throughout the paper complex manifolds are always assumed to be connected.

\section{Different notions of parabolicity}

\bigskip

\begin{definition}
A Stein manifold $\,\,X\,\,$ is called parabolic, in case it does not
possess a non-constant bounded above plurisubharmonic function.
\end{definition}

Thus, parabolicity of $\,\,X\,\,$ is equivalent to the following: if $%
u\left( z\right) \in psh\left( X\right) $ and $u\left( z\right) \leqslant C,$
then $u\left( z\right) \equiv const$ on $X$. It is convenient to describe
parabolicity in term of P-measures of pluripotential theory \cite{SA3}. We
will briefly recall this notion which can be defined for a general Stein
manifold $X.$ In the discussion below we will assume without loss a
generality that our Stein manifold $X$ is properly imbedded in $\mathbb{C}%
_{w}^{2n+1},n=\dim X$ , and $\sigma \left( z\right) $ is the restriction of $%
\ln \left\vert w\right\vert $ to $X$. Then $\sigma \left( z\right) \,\in
psh\left( X\right) \cap C^{\infty }\left( X\right) \,\,,$ $\left\{ {\sigma
\left( z\right) \leqslant C}\right\} \subset \subset X\,$\ $\forall C\in R$.
We further assume that $0\notin X\,$ and $\min \sigma \left( z\right) \,<0$.
We consider $\,\,{\sigma }$- balls $\,\,B_{R}=\left\{ {z\in X:\sigma \left(
z\right) <\ln R}\right\} {\text{ }}$ and as usual, define the class $\aleph
\left( {\overline{B}_{1},B_{R}}\right) \,\,,\,R>1$, of functions $u\left(
z\right) \in psh\left( {B_{R}}\right) $ such, that $u|_{B_{R}}\leqslant
0\,\,,\,\,u|_{\overline{B}_{1}}\leqslant -1.$ We put 
\begin{equation*}
\omega \left( {z\,,\,\overline{B}_{1}\,,\,B_{R}}\right) \,=\,\sup \left\{ {%
u\left( z\right) :u\in \aleph \left( {\overline{B}_{1},B_{R}}\right) }%
\right\} .
\end{equation*}%
The regularization $\omega ^{\ast }\left( {z\,,\,\overline{B}_{1}\,,\,B_{R}}%
\right) $ is called the P- measure of $\overline{B}_{1}$ with respect to $%
B_{R}$\cite{SA3}.

The P- measure $\omega ^{\ast }\left( {z\,,\,\overline{B}_{1}\,,\,B_{R}}%
\right) $ is plurisubharmonic in $B_{R}$, is equal to $-1$ on $\overline{B}%
_{1} $ and tends to 0 for $z\rightarrow \partial B_{R}.$ Moreover, it is
maximal, that is $\left( {dd^{c}\omega ^{\ast }}\right) ^{n}=0\,\,$ in $%
B_{R}\backslash \overline{B}_{1}$ and decreasing by $R$. We put $\omega
^{\ast }\left( {z,\overline{B}_{1}}\right) =\lim_{R\rightarrow \infty
}\omega ^{\ast }\left( {z,\overline{B}_{1},B_{R}}\right) .$

It follows that, 
\begin{equation*}
\omega ^{\ast }\left( {z,\overline{B}_{1}}\right) \in psh\left( X\right)
\,,\,\,\,\,-1\leqslant \omega ^{\ast }\left( {z,\overline{B}_{1}}\right)
\,\leqslant 0
\end{equation*}%
and $\left( dd^{c}\omega ^{\ast }\left( {z,\overline{B}_{1}}\right) \right)
^{n}=0$ in $B_{R}\backslash \overline{B}_{1}$.

In the construction of $\omega ^{\ast }\left( {z,\overline{B}_{1}}\right) $
we have used the exhaustion function $\sigma \left( z\right) ,$ however it
is not difficult to see that $\omega ^{\ast }\left( {z,\overline{B}_{1}}%
\right) $ depends only on $X$ and ${\overline{B}_{1}};$ not on the choice of
the exhaustion function. Indeed one can define the P-measure for any non
pluripolar compact $K\subset X,$ by selecting any sequence of domains 
\newline
$K\subset D_{j}\subset \subset D_{j+1}\subset \subset X\,,$\ $%
\,j=1,2,...,X=\bigcup_{j=1}^{\infty }{D_{j}}$ and employing the above
procedure with $B_{R}$ 's replaced with $D_{j}$ 's. It follows from the
definition, that the plurisubharmonic functions $\omega ^{\ast }\left( {z,K}%
\right) +1$ and $\omega ^{\ast }\left( {z,\overline{B}_{1}}\right) +1$ are
dominated by a constant multiple of each other. In particular $\omega ^{\ast
}\left( {z,K}\right) \equiv -1$ if and only if $\omega ^{\ast }\left( {z,%
\overline{B}_{1}}\right) \equiv -1.$ Hence the later property is an inner
property of $X.$ For further properties of P-measures we refer the reader to 
\cite{SA3}, \cite{SA5}, \cite{Z}, \cite{Z2}.

Vanishing \textit{\ of } $\omega ^{\ast }\left( {z,K}\right) +1$ on a
parabolic manifold not only imply the triviality of bounded holomorphic
functions but also give some information on the growth of unbounded
holomorphic functions. In fact on parabolic manifolds, a kind of "Hadamard
three domains theorem" with controlled exponents, is true. The precise
formulation of this characteristic, that will appear below, is an adaptation
of the property $\left( DN\right) $ of Vogt \cite{V1}, which was defined for
general Fr\'{e}chet spaces, to the Fr\'{e}chet spaces of analytic functions.
As usual we will denote by $O\left( X\right) $ the\ Fr\'{e}chet spaces of
analytic functions defined on a complex manifold $X$ with the topology of
uniform convergence on its compact subsets. The proposition\ we will give
below is due to Zaharyuta \cite{Z3} and it has been independently
rediscovered by several other authors \cite{AY2},\cite{V3}. We will include
a proof of this result for the convenience of the reader. \qquad

\begin{proposition}
\emph{\textbf{.}} The following are equivalent for a Stein manifold $X$

a) $X$is parabolic;

b) P-measures are trivial on $X$ \ i.e. $\omega ^{\ast }\left( {z,K}\right)
\equiv -1$ for every nonpluripolar compact $K\subseteq X;$ \ \ \ \ \ \ 

c\ ) For every nonpluripolar compact set $K_{0}\subset X$ \ and for every
compact set $K$\ of $X$ there is another compact set\ $L$\ containing $K$\
such that%
\begin{equation*}
\left\Vert f\right\Vert _{K}\leq \left\Vert f\right\Vert _{K_{0}}^{\frac{1}{2%
}} \left\Vert f\right\Vert _{L}^{\frac{1}{2}},\,\,\,\forall \text{ }f \in
O\left( X\right), \text{ \ \ \ \ \ \ \ \ \ \ \ }\left( DN\text{ condition of
Vogt}\right)
\end{equation*}
where $\left\Vert {\,\ast\,} \right\Vert _{H}$\ denotes the sup norm on $H.$
\end{proposition}

\begin{proof}
If $X$ is parabolic, then $\omega ^{\ast }\left( {z,\overline{B}_{1}}\right) 
$ being bounded and plurisubharmonic on $X$ reduces to $-1$.

Conversely, let $u\left( z\right) $ be an arbitrary bounded above $psh$
function on $X$. Let $u_{R\text{ }}=\sup_{B_{R}}u\left( z\right) \,,\,\infty
\geq R\geqslant 1$. If $u\left( z\right) \neq const$, then 
\begin{equation*}
\frac{{u\left( z\right) -u_{R}}}{{u_{R}-u_{1}}}\in \aleph \left( {\overline{B%
}_{1},B_{R}}\right) \,\,\,\,\text{and hence}\,\,\,\,\frac{{u\left( z\right)
-u_{R}}}{{u_{R}-u_{1}}}\leqslant \omega ^{\ast }\left( {z,\overline{B}%
_{1},B_{R}}\right) .
\end{equation*}%
It follows that 
\begin{equation*}
u\left( z\right) \leqslant -u_{1}\omega ^{\ast }\left( {z,\overline{B}%
_{1},B_{R}}\right) +u_{R}\left( {1+\omega ^{\ast }\left( {z,\overline{B}%
_{1},B_{R}}\right) }\right) \,,\,z\in B_{R}\,,\,\,\eqno(1)
\end{equation*}%
and as $R\rightarrow \infty $ this gives

\begin{equation*}
u\left( z\right) \leqslant -u_{1}\omega ^{\ast }\left( {z,\overline{B}_{1}}%
\right) +u_{\infty }\left( {1+\omega ^{\ast }\left( {z,\overline{B}_{1}}%
\right) }\right) \,,\,z\in X\,.\eqno(2)
\end{equation*}%
If $\omega ^{\ast }(z,\overline{B}_{1})\equiv -1$, then $u(z)\leq
u_{1},\,\,z\in X$, and by maximal principle we have $u\left( z\right)
=u_{1}\equiv const$, so that a) and b) are equivalent.

Now we look at the sup norms $\left\Vert {\,.\,}\right\Vert _{B_{m}}$ on the
sublevel balls $B_{m}\supset K$. Choose an increasing sequence of norms $%
\,\left\Vert {\,.\,}\right\Vert _{k}=\,\left\Vert {\,.\,}\right\Vert
_{m_{k}},\,\,k=0,1,\ldots ,\,\,m_{0}=1,$ that satisfy the condition c) with
the dominating norm $\left\Vert {\,.\,}\right\Vert _{0}$: 
\begin{equation*}
\,\,\left\Vert f\right\Vert _{k}\leqslant \,\left\Vert f\right\Vert _{{0}}^{%
\frac{1}{2}}\left\Vert f\right\Vert _{k+1}^{\frac{1}{2}}\,\,\,\,\forall f\in
O\left( X\right)
\end{equation*}

Iterating this inequality one gets: 
\begin{equation*}
\left\Vert f\right\Vert _{1}\leqslant \,\left\Vert f\right\Vert _{0}^{\frac{{%
2^{k-1}-1}}{{2^{k-1}}}}\,\left\Vert f\right\Vert _{k}^{\frac{1}{{2^{k-1}}}%
}\,\,\forall \,f\in \,O\left( X\right) .\eqno(3)
\end{equation*}

Denoting the sequence of domains by $D_{k}=B_{m_{k}}$ we consider the P-
measures $\,\omega ^{\ast }\left( {z,}\overline{D}_{0},D_{k+1}\right) $ , $\
\ k=1,2...$. Since these functions are continuous, by Bremermann's theorem
(see \cite{BR1}) for a fixed $k$ we can find analytic functions $%
f_{1},f_{2},....,f_{m}$ on $D_{k+1}$ and positive numbers $%
a_{1},a_{2},...,a_{m}$ such that 
\begin{equation*}
\omega ^{\ast }\left( z,\overline{D}_{0},D_{k+1}\right) \,+1-\varepsilon
\leqslant \max_{1\leqslant j\leqslant m}\left( {a_{j}\ln \left\vert {%
f_{j}\left( z\right) }\right\vert }\right) \,\leqslant \,\omega ^{\ast
}\left( z,\overline{D}_{0},D_{k+1}\right) \,+1
\end{equation*}%
pointwise on $\overline{D}_{k}$. We note that the compact $\overline{D}_{k}$
is polynomially convex in \newline
$\mathbb{C}^{2n+1}\supset X,$ \ so by Runge's theorem the functions $f_{j}$
can be uniformly approximated on $\overline{D}_{k}$ by functions $F\in 
\mathrm{O}\left( X\right) $. This in turn by (3) gives us the estimate $%
\omega ^{\ast }\left( z,\overline{D}_{0},D_{k+1}\right) \,+1\leqslant \frac{1%
}{{2^{k-1}}}+\varepsilon ,z\in D_{1}$. Now playing the same game with $D_{1}$
replaced by a given $D_{j}$ we see that $\omega ^{\ast }\left( {z,K,D_{k+1}}%
\right) \,$ converge uniformly to $-1$ on any compact subset of $X$, i.e. $%
\omega ^{\ast }\left( z,\overline{D}_{0}\right) \equiv -1$. Hence c) $%
\Rightarrow $ b).

Conversely, suppose that $\omega ^{\ast }\left( {z,K}\right) \equiv -1$ for
any $K\subset \subset X$. Fix a nonpluripolar compact $K_{0}\subset X$ and
fix an arbitrary compact set $K\subset X.$ Let $B_{k_{0}}\supset {\
K_{0}\cup K},\,\,k_{0}\in {\mathbb{N}}.$ Then, in view of Dini's theorem we
can choose $k$ so large that\newline
$\omega ^{\ast }\left( {z,K_{0},B_{k}}\right) \leqslant -1/2\,\,\,\,\text{for%
}\,\,z\in B_{k_{0}}.$ Since $\omega ^{\ast }\left( {z,K_{0},B_{k}}\right) $
is maximal on $B_{k}\backslash K_{0}$, then for arbitrary $\text{\thinspace }%
{f}\in O\left( {X}\right) ,f\neq 0,$ the inequality 
\begin{equation*}
\frac{{\ln \frac{{\left\vert {f\left( z\right) }\right\vert }}{{\left\Vert
f\right\Vert _{K_{0}}}}}}{{\ln \frac{{\left\Vert f\right\Vert _{B_{k}}}}{{%
\left\Vert f\right\Vert _{K_{0}}}}}}\,\,\leqslant \,\omega ^{\ast }\left( {%
z,K_{0},B_{k}}\right) \,+1\,,\,\,z\in B_{k},
\end{equation*}%
is valid. This in turn implies that 
\begin{equation*}
\ \left\Vert f\right\Vert _{K}\leqslant \left\Vert f\right\Vert
_{B_{k_{0}}}\leqslant \left\Vert f\right\Vert _{K_{0}}^{1/2}\left\Vert
f\right\Vert _{B_{k}}^{1/2},\,\,\,
\end{equation*}%
for all $f\in \mathrm{O}\left( X\right) .$ Hence $a)\Rightarrow c)$. This
finishes the proof of the proposition.
\end{proof}

\bigskip

\begin{definition}
A Stein manifold $X$ is called $S$ -parabolic, if there exit exhaustion
function $\rho \left( z\right) \in psh\left( X\right) $ that is maximal
outside a compact subset of $X.$ If in addition we can choose $\rho $ to be
continuous then we will say that $X$ is $S^{\ast }$-parabolic.
\end{definition}

In previous papers on parabolic manifolds (see for example \cite{GK},\cite%
{ST1}) authors usually required the condition of $C^{\infty }$ - smoothness
of $\rho $. Here we only distinguish the cases when the exhaustion function
is just $psh$ or continuous $psh$.

It is not difficult to see that $S$-parabolic manifolds are parabolic. In
fact, since the exhaustion function $\rho \left( z\right) $ of the
definition of $S$-parabolicity is maximal off some compact $K\subset \subset
X$, the balls $B_{r}=\left\{ {\rho \left( z\right) <\ln r}\right\}
,r\geqslant r_{0},$ contain $K$ for big enough $r_{0}$ and hence

\begin{equation*}
\omega ^{\ast }\left( {z,\overline{B}_{r_{0}},B_{R}}\right) =\max \left\{-1,%
\frac{\rho \left( z\right) -R}{R-r_{0}}\right\} .
\end{equation*}%
Consequently, 
\begin{equation*}
\omega ^{\ast }\left( {z,\overline{B}_{r_{0}}}\right) =\lim_{R\rightarrow
\infty }\omega ^{\ast }\left( {z,\overline{B}_{r_{0}},B_{R}}\right) \equiv
-1\,\,,\,\,z\in X.
\end{equation*}

For Stein manifolds of dimension one, the notions of $\ S$ -parabolicity, $%
S^{\ast }$-parabolicity, and parabolicity coincide. This is a consequence of
the existence of Evans-Selberg potentials ( subharmonic exhaustion functions
that are harmonic outside a given point) on a parabolic Riemann surfaces 
\cite{SAR}.

\bigskip

\begin{problem}
Do the notions of$\ \ S$ -parabolicity and $S^{\ast }$-parabolicity coincide
for Stein manifolds of arbitrary dimension?
\end{problem}

\begin{problem}
Do the notions of$\ \ $parabolicity and $S$ -parabolicity coincide for Stein
manifolds of arbitrary dimension?
\end{problem}

\section{\textbf{Spaces of Analytic Functions on Parabolic Manifolds}}

In this section we will relate the above discussed notions of parabolicity
of a Stein manifold $X$ with the linear topological structure of $O(X)$%
\textbf{, }the Fr\'{e}chet space of analytic functions on $X$\ with the
compact open topology. The first result which we will state is due to
Aytuna-Krone-Terzioglu and it characterizes parabolicity of a Stein manifold 
$X$ of dimension n, in terms of the similarity of the linear topological
structures of $O\left( X\right) $ and $O\left( {\mathbb{C}^{n}}\right) $.
For the proof, we refer the reader to \cite{AY3}.

\begin{theorem}
\emph{\textbf{\ }}For a Stein manifold $X$ of dimension n the following are
equivalent:

a) $X$ is parabolic;

b) $\mathrm{O}\left( X\right) $ is isomorphic as Fr\'{e}chet spaces to $%
\mathrm{O}\left( {\mathbb{C}^{n}}\right) $.
\end{theorem}

The correspondence that sends an entire function $f$ to its Taylor
coefficients $\left( x_{m}\right) _{m=0}^{\infty }$ ordered in the usual
way, establishes an isomorphism between $O\left( 
%TCIMACRO{\U{2102} }%
%BeginExpansion
\mathbb{C}
%EndExpansion
^{n}\right) $ and the infinite type power series space 
\begin{equation*}
\Lambda _{\infty }\left( \alpha _{m}\right) :=\left( x=\left( x_{m}\right)
_{m=0}^{\infty }:\left\vert x\right\vert _{k}:=\sum\limits_{m=0}^{\infty
}\left\vert x_{m}\right\vert e^{k\alpha _{m}}<\infty \text{ }\forall
k=1,2,....\right)
\end{equation*}%
with $\alpha _{m}=m^{\frac{1}{n}},$ $m=0,1,2,.$.

Recall that a \textit{graded Fr\'{e}chet space\ }is a tuple $\left( X,\left(
\left\vert \ast \right\vert _{s}\right) \right) ,$ where \ $X$ is a Fr\'{e}%
chet space and $\left( \left\vert \ast \right\vert _{s}\right)
_{s=1}^{\infty }$ is a fixed system of seminorms defining the topology of $%
X. $ Whenever we deal with $\Lambda _{\infty }\left( \alpha _{m}\right) ,$%
for an exponent sequence\textbf{\ }$\alpha _{m}\uparrow \infty $ (not
necessarily $\left( m^{\frac{1}{n}}\right) _{m}$) $,$we will tacitly assume
that we are dealing with a graded space and that the grading is given by the
norms defined in the above expression . We will need a definition from the
structure theory of Fr\'{e}chet spaces;

\begin{definition}
A continuous linear operator T between two graded Fr\'{e}chet spaces $\left(
X,\left( \left\vert \ast \right\vert _{k}\right) \right) $ and $\left(
Y,\left( \left\Vert \ast \right\Vert _{k}\right) \right) $ is tame in case: 
\begin{equation*}
\exists A>0\text{ }\forall \text{ }k\text{ }\exists C>0:\left\Vert T\left(
x\right) \right\Vert _{k}\leq C\left\vert x\right\vert _{k+A},\forall x\in X.
\end{equation*}%
Two graded Fr\'{e}chet spaces are called tamely isomorphic in case there is
a one to one tame linear operator from one onto the other whose inverse is
also tame.
\end{definition}

The graded space $\left( O\left( 
%TCIMACRO{\U{2102} }%
%BeginExpansion
\mathbb{C}
%EndExpansion
\right) ,\left\Vert \ast \right\Vert _{k}\right) $ where $\left\Vert \ast
\right\Vert _{k}$ is the sup norm on the disc with radius $e^{k}$, is tamely
isomorphic, under the correspondence described above, to the power series
space $\Lambda _{\infty }\left( n\right) $ in view of the Cauchy's
inequality. This observation motivates our next definition:

\begin{definition}
Let $X$ be a Stein manifold. The space $O\left( X\right) $ is said to be
tamely isomorphic to an infinite type power series space in case there is an
exhaustion of $M$ by connected holomorphically convex compact sets $\left(
K_{k}\right) _{k=1}^{\infty }$ with $K_{k}\subset K_{k+1}^{\circ }$ $%
,n=1,2,. $, such that the graded space $\left( O\left( X\right) ,\left(
\sup_{K_{k}}\left\vert \ast \right\vert \right) \right) $ is tamely
isomorphic to an infinite type power series space$.$
\end{definition}

The supremum norms are, in some sense, associated with the function theory
whereas the power series norms are associated with the structure theory of Fr%
\'{e}chet spaces, and tame equivalence gives one, a controlled equivalence
between these generating norm systems.

For a graded Fr\'{e}chet space, linear topological properties that ensure
tame equivalence to an infinite type power series space were studied by D.
Vogt and his school, in the context of structure theory of nuclear Fr\'{e}%
chet spaces \ \cite{V3}, \cite{POP}. In the proof of the theorem below, we
will make use of a specific result of Vogt in this direction. We recall this
theorem here for the benefit of the reader. '

\bigskip

\begin{definition}
A graded nuclear Fr\'{e}chet\textbf{\ }space\textbf{\ }$\left(
E,^{\left\vert \ast \right\vert _{k}}\right) $ is said to be a $\left(
DN\right) $ space in standard form in case, with suitable constants $%
C_{k}>0, $ 
\begin{equation*}
\left\vert \ast \right\vert _{k}^{2}\leq C_{k}\left\vert \ast \right\vert
_{k-1}\left\vert \ast \right\vert _{k+1}
\end{equation*}%
for all $k=1,2,...$
\end{definition}

\bigskip

\begin{definition}
A graded nuclear Fr\'{e}chet\textbf{\ }space\textbf{\ }$\left(
F,^{\left\vert \ast \right\vert _{k}}\right) $ is said to be a $\left(
\Omega \right) $ space in standard form in case, with suitable constants $%
D_{k}>0,$ 
\begin{equation*}
\left\vert \ast \right\vert _{k}^{\ast 2}\leq D_{k}\left\vert \ast
\right\vert _{k-1}^{\ast }\left\vert \ast \right\vert _{k+1}^{\ast }
\end{equation*}%
for all $k=1,2,...$where $\left\vert x^{\ast }\right\vert _{k}^{\ast }=\sup
\{\left\vert x^{\ast }\left( y\right) \right\vert :\left\vert y\right\vert
_{k}\leq 1\}$ $k=1,2,...,$denotes the dual "norms" on $F^{\ast }$
\end{definition}

\bigskip

\begin{theorem}
$\left( \text{cf \cite{V3}, Th 2.3}\right) $ Let $E$ be a nuclear $\left(
DN\right) $ space in standard form, and $F$ an $\left( \Omega \right) $
space in standard form. Suppose that there exits a tame surjection from $F$
onto $E$. Then $E$ is tamely isomorphic to an infinite type power series
space.
\end{theorem}

\bigskip

For two nonnegative real valued functions $\alpha $ and $\beta $ on a set $T$
we will use the notation $\alpha \left( t\right) \prec \beta \left( t\right) 
$ to mean $\exists $ $C>0$ such that $\alpha \left( t\right) \leq C\beta
\left( t\right) $ $\forall t\in T.$

We can now state the main result of this section:

\bigskip

\begin{theorem}
Let $X$ be a Stein manifold. The space of analytic functions on $X$, $%
O\left( X\right) ,$ is tamely isomorphic to an infinite type power series
space if and only if $X$ \ is $S^{\ast }-Parabolic.$
\end{theorem}

\begin{proof}
$\Rightarrow :$ \ Suppose that $\ O\left( X\right) $ is tamely isomorphic to
a power series space $\Lambda _{\infty }\left( \alpha _{m}\right) $\textbf{\ 
}with\textbf{\ }$\alpha _{m}\uparrow \infty $. Fix a tame isomorphism $%
T:\Lambda _{\infty }\left( \alpha _{m}\right) \rightarrow O\left( X\right) $%
. We also fix an exhaustion $\left\{ K_{k}\right\} _{k=0}^{\infty }$ of $X%
\mathbf{\ }$by holomorphically convex compact sets and an integer $B^{\prime
}$ such that for all $k$ 
\begin{equation*}
\left\Vert T\left( x\right) \right\Vert _{k}\prec \left\vert x\right\vert
_{k+B^{\prime }}\text{ and }\left\vert x\right\vert _{k}\prec \left\Vert
T\left( x\right) \right\Vert _{k+B^{\prime }}\text{ }\forall x\in \mathbf{%
\Lambda }_{\infty }\left( \alpha _{m}\right) ,
\end{equation*}%
where $\left\Vert \ast \right\Vert _{k}$ denotes the sup norm on $K_{k},$ $%
k=0,1,2..$. Let $e_{m}\circeq $ $T\left( \varepsilon _{m}\right) $ where ,as
usual, $\varepsilon _{m}=\left( 0,...,0,1,0,...\right) ,$ $m=1,2...$. Set 
\begin{equation*}
\rho \left( z\right) \circeq \lim \sup_{\xi \rightarrow z}\lim
\sup_{m\rightarrow \infty }\text{ }\frac{\log \left\vert e_{m}\left( \xi
\right) \right\vert }{\alpha _{m}}.
\end{equation*}%
Clearly $\rho $ is a plurisubharmonic function on $X$. If we set $D_{\tau
}\circeq \left\{ z:\rho \left( z\right) <\tau \right\} $ for $\tau \in 
%TCIMACRO{\U{211d} }%
%BeginExpansion
\mathbb{R}
%EndExpansion
$ ,we have: 
\begin{equation*}
K_{k}\subseteq D_{k+B}\text{ \ for large }n,\text{\ where }B=B^{\prime }+1.
\end{equation*}%
Now fix an arbitrary $z_{0}\in D_{\sigma }$ and choose, in view of Hartogs
lemma, a small $\epsilon >0$ such that $\left\vert e_{m}\left( z_{0}\right)
\right\vert \leq Ce^{\alpha _{m}\left( \sigma -\epsilon \right) }$ for all%
\textbf{\ }$m$ .\textbf{\ }For any $x=\sum x_{m}\varepsilon _{m}\in \Lambda
_{\infty }\left( \alpha _{m}\right) $ we have:%
\begin{equation*}
\left\vert T\left( x\right) \left( z_{0}\right) \right\vert \mathbf{\leq }%
\sum_{m}\left\vert x_{m}\right\vert \text{ }\left\vert e_{m}\left(
z_{0}\right) \right\vert \prec \sum_{m}\left\vert x_{m}\right\vert e^{\left(
\sigma -\epsilon \right) \alpha _{m}}\mathbf{\prec }\left\Vert T\left(
x\right) \right\Vert _{\left[ \left[ \sigma \right] \right] +1+B}\text{ }
\end{equation*}%
\textbf{\bigskip }Since T is onto and $K_{m}$'s are holomorphically convex,
we have that $z_{0}\in K_{_{\left[ \left[ \sigma \right] \right] +1+B}}$.
Hence $D_{\sigma }\subseteq K_{_{\left[ \left[ \sigma \right] \right] +1+B}}$%
. Combining this with our previous findings we get%
\begin{equation*}
\exists \text{ }d>0\text{ such that }\overline{D_{\sigma }}\subseteq
D_{\sigma +d}\text{ }\forall \alpha \text{ large}
\end{equation*}%
\bigskip Now fix a nice compact set $K,$ say $K=\overline{D}$ for some
domain, with the property that 
\begin{equation*}
\exists \text{\ }\beta >0\text{ such that }\left\vert x\right\vert _{\beta
}\prec \sup_{w\in K}\left\vert T\left( x\right) \mathbf{(}w\mathbf{)}%
\right\vert \text{ \ }\forall x\in \Lambda _{\infty }\left( \alpha
_{m}\right) .
\end{equation*}%
\bigskip We wish to show that 
\begin{equation*}
\Phi \left( z\right) \circeq \lim \sup_{\xi \rightarrow z}\left\{ \varphi
\left( \xi \right) :\varphi \in psh\left( X\right) ,\text{ }\varphi |K\leq
0,\varphi \leq \rho +C\text{ for some }C=C\left( \varphi \right) \right\}
\end{equation*}%
\bigskip defines a plurisubharmonic function on $X.$ To this end choose a%
\textbf{\ }$\varphi \in psh\left( X\right) $\textbf{\ }with\textbf{\ }$%
\varphi |K\leq 0$\textbf{\ }and $\varphi \leq \rho +C$\textbf{\ }for some%
\textbf{\ }$C=C\left( \varphi \right) >0$\textbf{.} By Bremermann's theorem\
( \cite{BR1}),we choose a representation 
\begin{equation*}
\varphi \left( z\right) =\lim \sup_{\xi \rightarrow z}\lim \sup_{j}\frac{%
\log \left\vert f_{j}\left( \xi \right) \right\vert }{c_{j}}
\end{equation*}%
of $\varphi $ on $X$ for some $f_{j}\mathbf{\in }O\left( X\right) \mathbf{\
\ }j=1,2\mathbf{....,}$ and positive real numbers\textbf{\ }$c_{j}\mathbf{%
\uparrow \infty ,}$ $j=1,2,3..$.. Using Hartogs lemma in a suitable
neighborhood of $K$ we get:%
\begin{equation*}
\mathbf{\forall }\text{ }\epsilon >0\text{ }\mathbf{\exists }\text{ }%
j_{0}:\left\vert f_{j}\left( x\right) \right\vert \text{ }\mathbf{\leq }%
e^{\epsilon c_{j}}\mathbf{,}\text{ }j\text{\textbf{\ }}\geq \text{ }j_{0},%
\text{ }z\in K
\end{equation*}%
\bigskip In particular if $y_{j}\circeq T^{-1}\left( f_{j}\right) $ we have: 
\begin{equation*}
\lim \sup_{j}\frac{\log \left\vert y_{j}\right\vert _{\beta }}{c_{j}}\leq 0.
\end{equation*}%
\bigskip

Taking into account the relation between $\varphi $\ and $\rho $\ and the
inclusion $\overline{D_{\sigma }}\subseteq D_{\sigma +d}$\ for large $\sigma 
$, we have, in view of Hartogs lemma\textbf{: }%
\begin{equation*}
\sup_{w\in D_{\sigma }}\left\vert f_{j}\left( w\right) \right\vert \prec
e^{\left( \sigma +d+C\right) c_{j}}\text{ \ \ \ }\forall j,\forall \sigma 
\text{ large}
\end{equation*}%
In particular for large $m,$ we have;%
\begin{equation*}
\left\vert y_{j}\right\vert _{m}\prec e^{\left( m+d+C+2B\right) c_{j}}\text{
,}\forall j.
\end{equation*}%
For\textbf{\ }any non negative number $t$, we define:%
\begin{equation*}
h\left( t\right) \circeq \lim \sup_{j}\frac{\log \left\vert y_{j}\right\vert
_{t}}{c_{j}}.
\end{equation*}%
This function is an increasing convex function on the positive real numbers.
Taking into account $h\left( \beta \right) \leqq 0$ and $h\left( m\right)
\leqq m+d+C+2B$ for large $m,$it follows, that 
\begin{equation*}
h\left( t\right) \leq \left( \frac{N+D}{N-\beta }\right) t-\left( \frac{N+D}{%
N-\beta }\right) \beta
\end{equation*}%
on the interval $\left[ \beta ,N\right] $ for every $N\geq \beta ,$where $%
D=m+d+C+2B$. \ Hence $h\left( t\right) \leq t-\beta $ for $t>>\beta .$

Going back, since 
\begin{equation*}
\sup_{w\in D_{\sigma }}\left\vert f_{j}\left( w\right) \right\vert \prec
\left\vert y_{j}\right\vert _{\sigma +2+2B}
\end{equation*}%
for $z$ with $\rho \left( z\right) =\sigma $, we see that,%
\begin{eqnarray*}
\varphi \left( z\right) &=&\lim \sup_{\xi \rightarrow z}\lim \sup_{n}\frac{%
\log \left\vert f_{n}\left( \xi \right) \right\vert }{c_{n}}\leq h\left(
\sigma +2+2B+d\right) \leq \sigma +2+2B+d-\beta \\
&=&\rho \left( z\right) +Q,\text{ where \ }Q=Q\left( B,d,\beta \right) \in 
%TCIMACRO{\U{211d} }%
%BeginExpansion
\mathbb{R}
%EndExpansion
^{+}.
\end{eqnarray*}

Hence 
\begin{equation*}
\Phi \left( z\right) \leq \rho \left( z\right) +Q
\end{equation*}%
and so $\Phi \in L_{loc}^{\infty }\left( X\right) $. In particular $\Phi $
is a plurisubharmonic function on $X$ and satisfies%
\begin{equation*}
\exists \text{ }C_{1\text{ }}>0\text{ and }C_{2}>0\text{ such that }\rho
\left( z\right) -C_{1}\leq \Phi \left( z\right) \leq \rho \left( z\right)
+C_{2}\text{ \ \ \ \ on }X.
\end{equation*}%
It follows that $\Phi $ is an exhaustion and as a free envelope, is maximal
outside a compact set \cite{BT}.

Observe also that the sublevel sets $\Omega _{r}\circeq \left\{ z:\Phi
\left( z\right) <r\right\} $ satisfy :%
\begin{equation*}
\exists \text{ }\kappa _{0}>0\text{ such that }\overline{\Omega _{r}}%
\subseteq \Omega _{r+\kappa _{0}}\text{ for }r\text{ large enough.}
\end{equation*}

Now fix a decreasing sequence$\left\{ u_{j}\right\} $ of continuous
plurisubharmonic functions on $X$ converging to $\Phi $. Fix a compact set $%
K $ and $\epsilon >0.$ Choose an $r$ so large that$\left( \frac{r+\text{ }%
\kappa _{0}-\frac{\epsilon }{2}}{r}-1\right) \max_{\xi \in K}\Phi \left( \xi
\right) \leq \frac{\epsilon }{2}.$ There exits an $j_{0}$ such that for $%
j\geq j_{0}$ on $\Omega _{r},$ $u_{j}\leq r+\kappa _{0}$ and $u_{j}|K\leq 
\frac{\epsilon }{2}.$ Hence on $\Omega _{r}:$%
\begin{equation*}
\frac{u_{j}-\frac{\epsilon }{2}}{r+\kappa _{0}-\frac{\epsilon }{2}}\leq
\omega ^{\ast }\left( K,\Omega _{r}\right) =\frac{1}{r}\Phi .
\end{equation*}%
where $\omega ^{\ast }$ is the corresponding P-measure (see section 2). It
follows that on $K,$ 
\begin{equation*}
0\leq u_{j}-\Phi \leq \left( \frac{r+\text{ }\kappa -\frac{\epsilon }{2}}{r}%
-1\right) \max_{\xi \in K}\Phi \left( \xi \right) +\frac{\epsilon }{2}\leq
\epsilon \text{ \ for \ }j\geq j_{0}\text{.}
\end{equation*}%
Hence the\ convergence is uniform on $K.$ So $\Phi $ is continuous.\bigskip

$\Leftarrow :$ \ In the proof of this implication will use the above
mentioned theorem of D. Vogt. However, first we wish bring to light a
particular $\Omega -$type condition for $O\left( X\right) $\ provided by a
given plurisubharmonic exhaustion function $\left( \text{see also \cite{Z}}%
\right) $. For this part of the argument, one does not need parabolicity. To
stress this point we will summarize our findings separately, in the below
Proposition .

Let $X$ be a Stein manifold and $\Phi $ $:X\rightarrow \lbrack -\infty
,\infty )$ a plurisubharmonic function that is an exhaustion. Set $D_{t}$ $%
=\left( x\text{ }|\Phi \left( x\right) <t\right) $ for $t\in 
%TCIMACRO{\U{211d} }%
%BeginExpansion
\mathbb{R}
%EndExpansion
.$ Choose an increasing function $\ell $ so that for each $t\in 
%TCIMACRO{\U{211d} }%
%BeginExpansion
\mathbb{R}
%EndExpansion
,$ $\overline{D_{t}}\subset D_{\ell \left( t\right) }.$ We fix a volume form
d$\mu $ on $X$ and using the notation of Lemma 1 \cite{AY1} , \ we let $%
d\varepsilon =cd\mu $\ where $c$\ is the strictly positive continuous
function that appears in Lemma 1 of \cite{AY1}\textbf{.} Set; 
\begin{equation*}
U_{t}=\left\{ f\in O\left( X\right) :\int_{D_{t}}\left\vert f\right\vert
^{2}d\varepsilon \leq 1\right\} .
\end{equation*}%
Fix positive numbers $s_{1},s_{2},s$ such that $\ell \left( 0\right)
<s_{1}\leq \ell \left( s_{1}\right) \leq s_{2}\leq \ell \left( s_{2}\right)
\leq s$ and $L\geq 0$ . Let 
\begin{equation*}
\text{\ }\Phi _{L}\left( z\right) \circeq \left( 
\begin{array}{c}
0\text{ \ \ if }\Phi \left( z\right) \leq 0 \\ 
\frac{L\Phi (z)}{s}\text{ \ otherwise\ \ \ \ \ \ }%
\end{array}%
\right. .
\end{equation*}%
Consider an analytic function $f\in U_{s_{2}}$ . \ Using Lemma 1 of \cite%
{AY1} , we choose a decomposition of $f$ on $W_{+}\cap $ $W_{-}$, $%
f=f_{+}-f_{-}$ , with $f_{\pm }\in O\left( W_{\pm }\right) $, $W_{+}=\left( 
\overline{D_{s_{1}}}\right) ^{c},$ $W_{-}=D_{s_{2}},$ and such that the
estimates 
\begin{equation*}
\int_{W_{\pm }}\left\vert f_{_{\pm }}\right\vert ^{2}e^{-\Phi
_{L}}d\varepsilon \leq K\int_{W_{+}\cap W_{-}}\left\vert f\right\vert
^{2}e^{-\Phi _{L}}d\mu
\end{equation*}%
hold with $K=K\left( X,s_{1,}s_{2,}s_{,}\Phi \right) >0.$ On the other hand,
since $f\in U_{s_{2}}$, 
\begin{equation*}
\int_{W_{+}\cap W_{-}}\left\vert f\right\vert ^{2}e^{-\Phi _{L}}d\mu \leq
C\int_{W_{+}\cap W_{-}}\left\vert f\right\vert ^{2}e^{-\Phi
_{L}}d\varepsilon \leq Ce^{-\frac{Ls_{1}}{s}}\text{ for some }C>0.
\end{equation*}%
Hence 
\begin{equation*}
\int_{W_{\pm }}\left\vert f_{_{\pm }}\right\vert ^{2}e^{-\Phi
_{L}}d\varepsilon \leq C_{1}e^{-\frac{Ls_{1}}{s}}\text{ \ \ \ \ for some }%
C_{1}>0.
\end{equation*}%
Now since\textbf{\ }$\Phi _{L}$\textbf{\ }is zero on\textbf{\ }$D_{0}$%
\textbf{\ }we have, 
\begin{equation*}
\int_{D_{0}}\left\vert f_{-}\right\vert ^{2}d\varepsilon
=\int_{D_{0}}\left\vert f_{-}\right\vert ^{2}e^{-\Phi _{L}}d\varepsilon \leq
\int_{W_{-}}\left\vert f_{-}\right\vert ^{2}e^{-\Phi _{L}}d\varepsilon \leq
C_{1}e^{-\frac{Ls_{1}}{s}}
\end{equation*}%
and 
\begin{equation*}
\int_{W_{-}}\left\vert f-f_{-}\right\vert ^{2}d\varepsilon \mathbf{\leq }%
C_{2}e^{\frac{L\left( s-s_{1}\right) }{s}}\mathbf{.}
\end{equation*}%
Set%
\begin{equation*}
G=\binom{f_{+}\text{ \ \ \ \ \ \ \ \ \ \ on }W_{+}}{f-f_{-}\text{ \ \ \ \ \
\ \ on }W_{-}}.
\end{equation*}%
Clearly\textbf{\ }$G\in O\left( X\right) ,$ and,%
\begin{eqnarray*}
\int_{D_{s}}\left\vert G\right\vert ^{2}d\varepsilon &\leq &\int_{D_{s}\cap
W_{+}}\left\vert G\right\vert ^{2}e^{-\Phi _{L}}e^{\Phi _{L}}d\varepsilon
+\int_{W_{-}}\left\vert G\right\vert ^{2}d\varepsilon \leq C_{3}\left( e^{%
\frac{L\left( s-s_{1}\right) }{s}}+e^{\frac{L\left( s-s_{1}\right) }{s}%
}\right) \\
&\leq &C_{4}e^{\frac{L\left( s-s_{1}\right) }{s}}.
\end{eqnarray*}%
Moreover%
\begin{equation*}
\int_{D_{0}}\left\vert G-f\right\vert ^{2}d\varepsilon
=\int_{D_{0}}\left\vert f_{-}\right\vert ^{2}d\varepsilon
=\int_{D_{0}}\left\vert f_{-}\right\vert ^{2}e^{-\Phi _{L}}d\varepsilon \leq
C_{1}e^{-\frac{Ls_{1}}{s}}.
\end{equation*}

Hence we obtain:%
\begin{equation*}
U_{s_{2}}\subseteq Ce^{-\frac{Ls_{1}}{s}}U_{0}+Ce^{\frac{L\left(
s-s_{1}\right) }{s}}U_{s}
\end{equation*}%
for some constant $C>0$ which does not depend upon $L.$

Set $t\circeq 1-\frac{s_{1}}{s},$ and $r=e^{L\left( 1-t\right) -\log C}.$
Varying the parameter $L$, a short computation yields%
\begin{equation*}
\exists \text{ }C>0\text{ such that: \ \ \ \ \ }U_{s_{2}}\subseteq \frac{1}{r%
}U_{0}+Cr^{\frac{t}{1-t}}U_{s}\text{ for all }r\in \left[ 1,\infty \right] .
\end{equation*}%
Since the above inclusion obviously holds for $0<r\leq 1,$ and writing the
value of $t$ we have:%
\begin{equation*}
\exists \text{ }D>0\text{ such that: \ \ \ \ \ }U_{s_{2}}\subseteq \frac{D}{r%
}U_{0}+\frac{r^{\frac{s}{s_{1}}}}{r}U_{s}\text{ for all }r\in \left(
0,\infty \right) .
\end{equation*}%
This is an $\Omega -$type condition introduced by Vogt and Wagner \cite{V2}
. In terms of the "dual norms" this condition can we expressed as ( {See} 
\cite{V2}):%
\begin{equation*}
\exists C>0\ \text{such that }\left\vert \left\vert x^{\ast }\right\vert
\right\vert _{s_{2}}^{\ast }\leq C\left( \left\vert \left\vert \left\vert
x^{\ast }\right\vert \right\vert \right\vert _{0}^{\ast }\right) ^{1-\frac{%
s_{1}}{s}}\left( \left\vert \left\vert \left\vert x^{\ast }\right\vert
\right\vert \right\vert _{s}^{\ast }\right) ^{\frac{s_{1}}{s}},\ \ \forall
x^{\ast }\in O\left( X\right) ^{\ast },
\end{equation*}

where $\left\vert \left\vert \left\vert x^{\ast }\right\vert \right\vert
\right\vert _{t}^{\ast }\circeq \sup \left\{ \left\vert x^{\ast }\left(
f\right) \right\vert :f\in O\left( X\right) ,\left\vert \left\vert
\left\vert f\right\vert \right\vert \right\vert _{t}\leq 1\right\} ,$ $\
x^{\ast }\in O\left( X\right) ^{\ast },$ $t\in R$ \ and$\left\vert
\left\vert \left\vert f\right\vert \right\vert \right\vert _{t}=\left(
\int_{D_{t}}\left\vert f\right\vert ^{2}d\varepsilon \right) ^{\frac{1}{2}}.$
We collect our findings, with the above notation, in:

\medskip \bigskip

\textbf{Proposition}: \ Let $X$\ be a Stein manifold and $\Phi $\ a
plurisubharmonic function on $X$ such that $D_{t}\circeq \left\{ z\text{ }%
|\Phi \left( z\right) <t\right\} \subset \subset X$\ $\ ,\forall t\in 
%TCIMACRO{\U{211d} }%
%BeginExpansion
\mathbb{R}
%EndExpansion
.$\ \ If we have 
\begin{equation*}
\overline{D_{s_{0}}}\subseteq D_{s_{1}}\subseteq \overline{D_{s_{1}}}%
\subseteq D_{s_{2}}\subseteq \overline{D_{s_{2}}}\subseteq D_{s}
\end{equation*}%
for some indexes $\ s_{0}<$\ $s_{1}<s_{2}<s,$\ then the Fr\'{e}chet space $%
O\left( X\right) ,$\ with the norms defined above, satisfies the following $%
\Omega -$\ condition:%
\begin{equation*}
\exists C>0\ :\left\vert \left\vert \left\vert x^{\ast }\right\vert
\right\vert \right\vert _{s_{2}}^{\ast }\leq C\left( \left\vert \left\vert
\left\vert x^{\ast }\right\vert \right\vert \right\vert _{s_{0}}^{\ast
}\right) ^{\frac{s-s_{1}}{s-s_{0}}}\left( \left\vert \left\vert \left\vert
x^{\ast }\right\vert \right\vert \right\vert _{s}^{\ast }\right) ^{\frac{%
s_{1}-s_{0}}{s-s_{0}}},\ \ \forall x^{\ast }\in O\left( X\right) ^{\ast }
\end{equation*}

\medskip \bigskip

Now we return to the proof of the theorem. Lets fix a continuous proper
plurisubharmonic function $\Phi $ on $X$ that is maximal outside a compact
set. We can arrange things so that $\Phi $ is maximal outside a compact
subset of $D_{0}$, where as usual $D_{t}=\left\{ x\text{ }|\text{ }\Phi
\left( x\right) <t\right\} .$ Since $\Phi $ is continuous, for a given $k$,
by taking $s_{0}=k-1-\frac{1}{k-1},$ $s_{2}=k-\frac{1}{k}$ , $s=k+1-\frac{1}{%
k+1}$ and choosing $s_{1}$, $s_{0}<s_{1}<s_{2},$ so that $\frac{s-s_{1}}{%
s-s_{0}}\leq \frac{1}{2}$, the above proposition gives%
\begin{equation*}
\forall k\geq 2\text{ }\exists \text{ }D_{k}>0\ :\left\vert \left\vert
\left\vert x^{\ast }\right\vert \right\vert \right\vert _{k-\frac{1}{k}%
}^{\ast }\leq C_{k}\left( \left\vert \left\vert \left\vert x^{\ast
}\right\vert \right\vert \right\vert _{k-1-\frac{1}{k-1}}^{\ast }\right) ^{%
\frac{1}{2}}\left( \left\vert \left\vert \left\vert x^{\ast }\right\vert
\right\vert \right\vert _{k+1-\frac{1}{k+1}}^{\ast }\right) ^{\frac{1}{2}},\
\ \forall x^{\ast }\in O\left( X\right) ^{\ast }
\end{equation*}%
Hence\ $O\left( X\right) $, with the grading $\left\vert \left\vert
\left\vert f\right\vert \right\vert \right\vert _{k}=\left( \int_{D_{k-\frac{%
1}{k}}}\left\vert f\right\vert ^{2}d\varepsilon \right) ^{\frac{1}{2}}$, $%
k=1,2....$ is an $\Omega -$space in \textit{standard} form. On the other
hand the grading\textbf{\ }$\left\Vert f\right\Vert _{k}=\sup_{z\in
D_{k}}\left\vert f\left( z\right) \right\vert $\textbf{\ , }$k=0,1,2...$,on $%
O\left( X\right) $, satisfies 
\begin{equation*}
\left\Vert f\right\Vert _{k}^{2}\leq \left\Vert f\right\Vert
_{k+1}\left\Vert f\right\Vert _{k-1}\text{ }
\end{equation*}%
In fact for a non constant analytic function $f$\ on $X$\ ,and fixed $%
k\varepsilon N$ ,the plurisubharmonic function 
\begin{equation*}
\rho \left( z\right) =2\frac{\log \left( \frac{\left\vert f\left( z\right)
\right\vert }{\left\Vert f\right\Vert _{k+1}}\right) }{\log \left( \frac{%
\left\vert f\left( z\right) \right\vert _{k+1}}{\left\Vert f\right\Vert
_{k-1}}\right) }
\end{equation*}%
is dominated by the maximal function $\Phi -\left( k+1\right) $\ on the
boundary of the region $\left( z:k-1<\Phi \left( z\right) <k+1\right) $\ \
and hence is dominated by it on the whole region. Rewriting this domination
on the level set $\Phi =k$\ yields the desired inequality.

Hence, $O\left( X\right) $\ with the grading $\left\Vert f\right\Vert
_{k}=\sup_{D_{k}}\left\vert f\right\vert $, $k=0,1,2..$\ is a $\left( 
\mathbf{\ }DN\right) -$\ space in standard form\textbf{.}

Moreover for every $k=1,2...$, there is a $K_{k}>0,$such that $\left\vert
\left\vert \left\vert f\right\vert \right\vert \right\vert _{k}\leq
K_{k}\left\Vert f\right\Vert _{k}$ and $\left\Vert f\right\Vert _{k}\leq
K_{k}\left\vert \left\vert \left\vert f\right\vert \right\vert \right\vert
_{k+2}$ . Now all the conditions of Vogt's theorem mentioned above, are
satisfied with identity as the required surjection. It follows that $O\left(
X\right) $ is tamely isomorphic to an infinite type power space. This
finishes the proof of the theorem.\bigskip
\end{proof}

\bigskip

\bigskip

The theorem above associates to every special plurisubharmonic continuous
exhaustion function $\Phi $\ on a $S^{\ast }-parabolic$\ \ Stein manifold\ $%
X,$\ an exponent sequence $\left( \alpha _{m}\right) _{m}$\ such that the
spaces \ $\left( O\left( X\right) ,^{\left\Vert \ast \right\Vert
_{k}}\right) $\ with grading coming from the sup norms on the level sets of $%
\Phi ,$\ and $\Lambda _{\infty }\left( \alpha _{m}\right) $\ are tamely
isomorphic. It might be of interest to examine the exponent sequences $%
\left( \alpha _{m}\right) _{m=0}^{\infty }$\ obtained in this way and see
how they depend upon the special exhaustion function $\Phi $\ $.$\ 

To this end let $X$\ be a Stein manifold with a continuous plurisubharmonic
exhaustion function $\Phi $\ that is maximal off a compact set that lies in
the interior of $K_{0}=\left( z:\Phi \left( z\right) \leq 0\right) .$\ We
will choose a hilbertian grading $\left( \left\vert \left\vert \ast
\right\vert \right\vert _{k}^{\wedge }\right) _{k}$\ of $O\left( X\right) $\
so that the Hilbert spaces $H_{k}\circeq \overline{\left( O\left( X\right)
,^{\left\vert \left\vert \ast \right\vert \right\vert _{k}^{\wedge }}\right) 
}$\ $_{k=0}^{\infty }$\ satisfy the continuous inclusions;%
\begin{equation*}
H_{k}\hookrightarrow O\left( D_{k}\right) \hookrightarrow A\left(
K_{0}\right) \hookrightarrow H_{0}\text{ \ \ \ }\forall k=1,2..
\end{equation*}%
where $D_{k}=\left( z:\Phi \left( z\right) <k\right) ,$and $A\left(
K_{0}\right) $\ is the germs of analytic functions on $K_{0}$\ with the
inductive topology. Moreover we also require that:

$\mathbf{\bigskip }$

$a)$ The tuple $\left( H_{0},H_{k}\right) $ is admissible for the pair $%
\left( K_{0},D_{k}\right) $in the sense of Zaharyuta \cite{Z3}$,\forall
k\epsilon \mathbb{N}$,

$\bigskip $

$b)$ The theorem above is valid i.e. there is an infinite type power series
space $\Lambda _{\infty }$ $\left( \alpha \right) $ so that $\ \left(
O\left( X\right) ,^{^{\left\vert \left\vert \ast \right\vert \right\vert
_{k}^{\wedge }}}\right) _{k=0}^{\infty }$ is tamely isomorphic to $\Lambda
_{\infty }$ $\left( \alpha \right) $ .

\bigskip

We will only use a special property of admissible pairs, so we will just
refer the reader to \cite{Z4} for the definition, construction and a
detailed discussion of this notion. However we should mention that in our
case we can take $H_{0}$\ to be the closure in $L^{2}\left( X,dd^{c}\max
\left( 0,\Phi \right) \right) ^{n})$, of the space of analytic functions
defined near $K_{0}$ , and $H_{k}$ to be $O\left( D_{k}\right) \cap
L^{2}\left( d\varepsilon \right) $ where $d\varepsilon $\ is the measure
that appears in the proof of Theorem 2 ( \cite{Z4} \cite{AY1}) and the
existence of an infinite type power series space\ satisfying the required
property for this choice of generating norms follows from the proof the
theorem given above. In what follows, we will denote the corresponding
graded space by $\left( O\left( X\right) ,^{\Phi }\right) .$

Since $O\left( X\right) ,$ for a parabolic\ Stein manifold\ $X$ of dimension 
$n,$ is isomorphic to $\Lambda _{\infty }\left( m^{\frac{1}{n}}\right) ,$
regardless of the special exhaustion function we have:%
\begin{equation*}
\exists C>0:\frac{1}{C}\leq \lim \inf_{m}\ \frac{\alpha _{m}}{m^{\frac{1}{n}}%
}\leq \lim \sup_{m}\ \frac{\alpha _{m}}{m^{\frac{1}{n}}}\leq C
\end{equation*}%
for all exponent sequences $\left( \alpha _{m}\right) _{m=0}^{\infty }$\
such that $O\left( X\right) $\ and $\Lambda _{\infty }$\ $\left( \alpha
\right) $\ are isomorphic.

To proceed further we need the notion of a Kolmogorov diameter. For a vector
space $L$, let us denote the collection of all subspaces of $Y\subset L$
with dim$Y\leq m,$ by $L_{m},m=1,2...$

\begin{definition}
Let $\left( \mathbb{X},^{\left\vert * \right\vert _{k}}\right) $ be a graded
Fr\'{e}chet space with an increasing sequence of seminorms. Let $%
U_{i}=\left( x\epsilon \mathbb{X}:\left\vert x\right\vert _{i}\leq 1\right)
,i=1,2...$. The $m^{th}$ diameter of $U_{i}$ with respect to $U_{j}$, $i<j,$
is defined by 
\begin{equation*}
d_{m}\left( U_{i},U_{j}\right) \circeq \inf \left( \lambda >0:\exists \text{ 
}Y\epsilon \text{ }\mathbb{X}_{k}\text{ such that }U_{i}\subseteq \lambda
U_{j}+Y\right) .
\end{equation*}
\end{definition}

\bigskip

Now fix a $S^{\ast }-parabolic$ \ Stein manifold\ $X$ and suppose that $%
\left( O\left( X\right) ,^{\Phi }\right) $ and $\Lambda _{\infty }\left(
\alpha _{m}\right) $ are tamely isomorphic under an isomorphism $T.$ In
particular\textbf{\ }there exits an $A>0$ such that, 
\begin{equation*}
\forall k\text{ }\exists C>0:\left\vert \left\vert T\left( x\right)
\right\vert \right\vert _{k}^{\wedge }\leq C\left\vert x\right\vert _{k+A}%
\text{ \ and }C\left\vert \left\vert T\left( x\right) \right\vert
\right\vert _{k+A}^{\wedge }\geq \left\vert x\right\vert _{k}\text{ ,\ }%
\forall x\epsilon \Lambda _{\infty }\left( \alpha _{m}\right) .
\end{equation*}

We will denote by $U_{i}$ and $V_{i}$ the unit balls corresponding to the $%
i^{th}$ norms of $\left( O\left( X\right) ,^{\Phi }\right) $ and $\Lambda
_{\infty }\left( \alpha _{m}\right) $ respectively.

Fix a $k>>l$ large and suppose 
\begin{equation*}
U_{k}\subseteq \lambda U_{l}+L,
\end{equation*}%
for some $\lambda >0$ and $L$ some m-dimensional subspace of $O\left(
X\right) .$ Applying $T^{-1}$ to both sides and using the tame continuity
estimates we have:

\begin{equation*}
\frac{1}{C}V_{k+A}\subseteq T^{-1}\left( U_{k}\right) \subseteq \lambda
T^{-1}\left( U_{l}\right) +L^{\prime }\subseteq \lambda CV_{l-A}+L^{\prime },%
\text{ \ \ \ \ \ }L^{\prime }\circeq T^{-1}\left( L\right) .
\end{equation*}

Hence 
\begin{equation*}
d_{m}\left( V_{k+A},V_{l-A}\right) \leq Cd_{m}\left( U_{k},U_{l}\right)
\end{equation*}%
for all m, where the constant depends only on indices $k$ and $l$.

Arguing in a similar fashion, we also have 
\begin{equation*}
d_{m}\left( U_{k+A},U_{l-A}\right) \leq Cd_{m}\left( V_{k},V_{l}\right) ,%
\text{ }\forall m
\end{equation*}

\bigskip

It is a standard fact that $d_{m}\left( V_{k},V_{l}\right) =e^{\left(
l-k\right) \alpha _{m}}$ for $\ k>>l$ (\cite{DU1}). On the other hand our
requirement of admissibility of the norms $\left( \left\vert \left\vert \ast
\right\vert \right\vert _{k}^{\wedge }\right) _{k}$ gives , in view of a
result of Nivoche-Poletsky-Zaharyuta (Theorem 5 of \cite{Z4} , see also, 
\cite{Ni}) the asymptotics%
\begin{equation*}
\lim_{m}\frac{-\ln d_{m}\left( U_{k},U_{l}\right) }{m^{\frac{1}{n}}}=\frac{%
2\pi \left( n!\right) ^{\frac{1}{n}}}{\left( C\left( \overline{D_{l}}%
,D_{k}\right) \right) ^{\frac{1}{n}}}\text{ \ \ \ }\forall k>>l
\end{equation*}%
where $D_{s}=\left( z:\Phi \left( z\right) <s\right) $\textbf{\ is} as
above, and $C\left( \overline{D_{l}},D_{k}\right) $ is the Bedford-Taylor
capacity of the condenser $\left( \overline{D_{l}},D_{k}\right) $\cite{BT}.

Putting all these things together we have:

\begin{eqnarray*}
\lim \inf_{m}\frac{\alpha _{m}}{m^{\frac{1}{n}}} &\geq &\lim_{m}\left[ \frac{%
-\ln d_{m}\left( U_{k},U_{l}\right) }{m^{\frac{1}{n}}}\left( \frac{-\ln C}{%
\left( k-l+2A\right) \left( -\ln d_{m}\left( U_{k+A},U_{l-A}\right) \right) }%
+\frac{1}{\left( k-l+2A\right) }\right) \right] = \\
&&\frac{2\pi \left( n!\right) ^{\frac{1}{n}}}{\left( C\left( \overline{D_{l}}%
,D_{k}\right) \right) ^{\frac{1}{n}}}.\frac{1}{\left( k-l+A\right) }.
\end{eqnarray*}

\bigskip 
\begin{eqnarray*}
\lim \sup_{m}\frac{\alpha _{m}}{m^{\frac{1}{n}}} &\leq &\lim_{m}\left[ \frac{%
-\ln d_{m}\left( U_{k+A},U_{l-A}\right) }{m^{\frac{1}{n}}}\left( \frac{\ln C%
}{\left( k-l\right) \left( -\ln d_{m}\left( U_{k+A},U_{l-A}\right) \right) }+%
\frac{1}{\left( k-l\right) }\right) \right] = \\
&&\frac{2\pi \left( n!\right) ^{\frac{1}{n}}}{\left( C\left( \overline{%
D_{l-A}},D_{k+A}\right) \right) ^{\frac{1}{n}}}.\frac{1}{\left( k-l\right) }.
\end{eqnarray*}

On the other hand, since $\Phi $ is maximal off a compact set we can use the
function 
\begin{equation*}
\rho \left( z\right) =\frac{\Phi -r}{r-s}
\end{equation*}%
to compute the capacity of the condenser $\left( \overline{D_{s}}%
,D_{r}\right) $ for $r\gg s$ large enough . To be precise, in our case we
get \cite{BT}: 
\begin{equation*}
C\left( \overline{D_{s}},D_{r}\right) =\frac{1}{\left( r-s\right) ^{n}}%
\int_{X}\left( dd^{c}\Phi \right) ^{n}.
\end{equation*}

Taking this into account, we obtain: 
\begin{equation*}
\lim_{m}\frac{\alpha _{m}}{m^{\frac{1}{n}}}=2\pi \left( n!\right) ^{\frac{1}{%
n}}\left( \int_{X}\left( dd^{c}\Phi \right) ^{n}\right) ^{-\frac{1}{n}}.
\end{equation*}

We collect our findings in the proposition below. As usual $\left\vert
\left\vert * \right\vert \right\vert _{K\text{ \ \ }}$denote the sup norm on
a given compact set $K.$

\begin{proposition}
Let $X$ be a $S^{\ast }-parabolic$ Stein manifold of dimension $n$. Fix a
plurisubharmonic exhaustion function\ $\Phi $ on $X$ that is maximal outside
a compact set. Then the exponent sequence $\left( \alpha _{m}\right) _{n}$
of the infinite type power series space associated to $X$ by Theorem 2 above
satisfies: 
\begin{equation*}
\lim_{m}\frac{\alpha _{m}}{m^{\frac{1}{n}}}=2\pi \left( n!\right) ^{\frac{1}{%
n}}\left( \int_{X}\left( dd^{c}\Phi \right) ^{n}.\right) ^{-\frac{1}{n}}
\end{equation*}
\end{proposition}

\bigskip

Note that one can construct a new plurisubharmonic exhaustion function that
is again maximal off a compact set and with a prescribed positive right hand
side value in the equation above, by simply multiplying the given exhaustion
function with a positive constant. In particular we have\textbf{\ : }

\begin{corollary}
A Stein manifold $X$ of dimension $n$ is $S^{\ast }-parabolic$ if and only
if there exits an exhaustion of $X$ by connected , holomorphically convex
compact sets $\left( K_{k}\right) _{k=1}^{\infty },$ $K_{k}\subset $ $\left(
K_{k+1}\right) ^{\circ },$ $k=1,2..,$ such that the graded spaces $\left(
O\left( X\right) ,^{\left\vert \left\vert \ast \right\vert \right\vert
_{K_{k}}}\right) $ is tamely isomorphic to $\left( O\left( \mathbb{C}%
^{n}\right) ,^{\left\vert \left\vert \ast \right\vert \right\vert _{\Delta
_{k}}}\right) $ , where $\Delta _{k}$ is the polydisc in $\mathbb{C}^{n}$
with radius $k.$
\end{corollary}

\bigskip

\section{\textbf{\ Some classes of parabolic manifolds}}

An immediate class of parabolic manifolds can be obtained by considering
Stein manifolds that admit a proper analytic surjection onto some $\mathbb{C}%
^{n}.$ Affine algebraic manifolds belong to this class. Moreover such
manifolds are $S^{\ast }$-parabolic \cite{ST2}.

In this section we will look at some ways of generating parabolic manifolds
and give some nontrivial examples. \bigskip

\textbf{1. Complements of analytic multifunctions.}

\bigskip $_{{}}$

Let $A\subset \mathbb{C}^{n}$ be a closed pluripolar set whose complement is
pseudoconvex. Such sets are called "analytic multifunctions" by some
authors. They are studied extensively by various authors and are extremely
important in approximation theory, in the theory of analytic continuation
and in the description of polynomial convex hulls (see \cite{AW}, \cite{BR}, 
\cite{N}, \cite{O}, \cite{SA7}, \cite{SLOD}, \cite{SLOD2}, \cite{Y} and
others). These sets are removable for the class of bounded plurisubharmonic
functions defined on their complements. Hence their complements are
parabolic Stein manifolds. We would like to restate Problem 2 given above in
this setting since we hope that it will be more tractable.

\begin{problem}
Let $A$ be as above. Is $X=\mathbb{C}^{n}\backslash A$, $\ S-parabolic$%
?\bigskip
\end{problem}

In classical case, $n=1$, every closed polar set $A\subset \mathbb{C}$ is
analytic multifunction. As is well-known, if $K\subset \subset {\mathbb{C}}$
is a closed polar set, then there exist a subharmonic in $\,\,{\mathbb{C}}%
\,\,$ and harmonic in $\mathbb{C}\backslash K$ function $\,\,u(z),\,\,$such
that $\,\,u|_{K}\equiv -\infty \,\,$ and $\,\,u(z)-\ln \left\vert
z\right\vert \rightarrow 0\,\,$ as $\,\,z\rightarrow \infty $. One can use
such functions to construct a special exhaustion function on $\mathbb{C}%
\backslash A.\,$ To this end fix a $z_{0}\notin K\circeq A\cup \{\infty \}$
an arbitrary point. Then there exist $u(z)\in sh(\overline{\mathbb{C}}%
\backslash \{z_{0}\})\cap har(\{\overline{\mathbb{C}}\backslash
K\}\backslash \{z_{0}\}):\left. u\right\vert _{K}\equiv -\infty $ and $%
u(z)\rightarrow +\infty $ as $z\rightarrow z_{0}$. Therefore, $\rho
(z)=-u(z) $ is exhaustion for $X=\mathbb{C}\backslash A$, with one singular
point $z_{0}$.

\bigskip

On the other hand if $A=\{p(z)=0\}\subset \mathbb{C}^{n}$ is an algebraic
set, then it is easy to see that the function 
\begin{equation*}
\rho (z)\circeq -\frac{1}{{\deg p}}\ln \left\vert p\right\vert +2\ln
\left\vert z\right\vert
\end{equation*}%
is \ a special exhaustion function for $\mathbb{C}^{n}\backslash A$ \cite%
{ZE1}.

\bigskip

\bigskip

\begin{theorem}
Let $A=\{F(^{\prime }z,z_{n})=z_{n}^{k}+f_{1}(^{\prime
}z)z_{n}^{k-1}+...+f_{k}(^{\prime }z)=0\}$ be a Weierstrass polynomial
(algebraiodal) set in $\mathbb{C}^{n}$, where $f_{j}\in O(\mathbb{C}^{n-1})$
are entire functions, $j=1,2,...,k$, $k\geqslant 1$. Then $X=\mathbb{C}%
^{n}\backslash A$ \ is S*-parabolic.
\end{theorem}

\begin{proof}
We put 
\begin{equation*}
\rho (z)=-\ln \left\vert {F(z)}\right\vert +\ln (\left\vert {^{\prime }z}%
\right\vert ^{2}+\left\vert {F(z)-1}\right\vert ^{2}).\eqno(11)
\end{equation*}%
Then $\rho (z)=-\infty $ precisely on the finite set $Q=\left\{ {^{\prime
}z=0,F\left( {^{\prime }0,z_{n}}\right) =1}\right\} .$ Moreover, $\rho $ is
maximal, $\left( {dd^{c}\rho }\right) ^{n}=0$ and continuous outside of $%
A\cup Q,$ because $-\ln \left\vert {F(z)}\right\vert $ is pluriharmonic and $%
\ln (\left\vert {^{\prime }z}\right\vert ^{2}+\left\vert {F(z)-1}\right\vert
^{2})$ is maximal, since for any holomorphic vector-function $%
f=(f_{1},f_{2},...,f_{n}):f\neq 0$ the function $\ln \left\Vert {f}%
\right\Vert ^{2}$ is a maximal $psh$ function outside the zero set of $f$.

We will show, that $\rho \left( z\right) $ is exhaustion on $X=\mathbb{C}%
^{n}/A$, i.e. 
\begin{equation*}
\left\{ {\rho \left( z\right) <R}\right\} \subset \subset X\,\,\text{%
for\thinspace \thinspace \thinspace every}\,\,R\in {\mathbb{R}}.\eqno(12)
\end{equation*}%
If $F(z)=0$, then $\rho (z)=+\infty +\ln (\left\vert {^{\prime }z}%
\right\vert ^{2}+1)=+\infty $, so that $\rho |_{A}=+\infty $. The condition
(12) is clear, if all $f_{j},$ $j=0,1,...,k,$ are constant and so we assume
that, at least one of them is not constant. Then $M_{R}=\max_{\left\vert {%
^{\prime }z}\right\vert \leqslant R}\{\left\vert {f_{1}(^{\prime }z)}%
\right\vert ,...,\left\vert {f_{k}(^{\prime }z)}\right\vert \}\rightarrow
\infty $. For $\left\vert {^{\prime }z}\right\vert =R\geqslant 1$ and $%
\left\vert {z_{n}}\right\vert \leqslant M_{R}^{2}$ we have 
\begin{equation*}
\rho (z)=\ln \frac{{\left\vert {^{\prime }z}\right\vert ^{2}+\left\vert {%
F(z)-1}\right\vert ^{2}}}{{\left\vert {F(z)}\right\vert }}\geqslant \ln 
\frac{{\left\vert {^{\prime }z}\right\vert ^{2}+\left\vert {F(z)-1}%
\right\vert ^{2}}}{{1+\left\vert {F(z)-1}\right\vert }}\geqslant \ln \frac{{%
\left\vert {^{\prime }z}\right\vert ^{2}+\left\vert {F(z)-1}\right\vert ^{2}}%
}{{\left\vert {^{\prime }z}\right\vert +\left\vert {F(z)-1}\right\vert }}%
\geqslant
\end{equation*}%
\begin{equation*}
\geqslant \ln \frac{{\left\vert {^{\prime }z}\right\vert +\left\vert {F(z)-1}%
\right\vert }}{2}\geqslant \ln \frac{R}{2}.
\end{equation*}%
On the other hand on $\left\vert {^{\prime }z}\right\vert \leqslant R\,\,\,%
\text{and}\,\,\left\vert {z_{n}}\right\vert =M_{R}^{2}$ \ we have:%
\begin{equation*}
\rho \left( z\right) =\ln \frac{{\left\vert {^{\prime }z}\right\vert
^{2}+\left\vert {F\left( z\right) -1}\right\vert ^{2}}}{{\left\vert {F\left(
z\right) }\right\vert }}\geq \ln \frac{{%
(M_{R}^{2k}-M_{R}M_{R}^{2k-2}-...-M_{R}-1)^{2}}}{{%
M_{R}^{2k}+M_{R}M_{R}^{2k-2}+...+M_{R}}}=
\end{equation*}%
\begin{equation*}
=\ln M_{R}^{2k}\left( {1+\alpha _{k}}\right) ,
\end{equation*}%
where $\alpha _{k}\rightarrow 0\,\,\text{for}\,\,R\rightarrow \infty .$ It
follows that $\,\rho |_{\partial U_{R}}\rightarrow +\infty \,\,\,\text{for}%
\,\,\,R\rightarrow \infty $, where $U_{R}=\left\{ {\left( {\left\vert {%
^{\prime }z}\right\vert \leqslant R,\left\vert {z_{n}}\right\vert \leqslant
M_{R}^{2}}\right) }\right\} $.

Let us now consider the level set $D_{C}=\left\{ {\rho \left( z\right) <C}%
\right\} \,,C-$ constant. It is an open set and it contains the pole set Q.
If $\,\,R\,\,$is so big, that $U_{R}\supset Q\,\,\,\text{and}$ $\min \left\{ 
{\ln \frac{R}{2},\ln M_{R}^{2k}\left( {1+\alpha _{R}}\right) }\right\}
\geqslant C,\,\,\text{then}\,\,D_{C}\subset \subset U_{R}$ , since $D_{C}$
has no any component outside $\,\,U_{R}\,\,$ because of maximality of $%
\,\,\rho \,\,$\ on $X\backslash U_{R}.$ This completes the proof that $%
\,\,\rho \,\,$ is an exhaustion function.\textbf{\ }
\end{proof}

\bigskip

\begin{corollary}
The complement , $\mathbb{C}^{n}/\Gamma $ , of the graph $\ \Gamma =\left\{
\left( {^{\prime }z,z_{n}}\right) \varepsilon \mathbb{C}^{n}\text{: }%
z_{n}=f\left( {^{\prime }z}\right) \right\} $ of an entire function $f$ is $%
S^{\ast }$-parabolic.
\end{corollary}

\bigskip

\textbf{2. Manifolds, which admit an exhaustion function with small $%
(dd^{c})^{n}\text{mass.}$}

\bigskip

Demailly \cite{D} considered manifolds X \ which admit a continuous
plurisubharmonic exhaustion function $\,\,\varphi ,\,\,$with the property
that, 
\begin{equation*}
\lim_{r\rightarrow \infty }\frac{{\int_{B_{r}}{\left( {dd^{c}\varphi }%
\right) ^{n}}}}{{\ln r}}\,=\,0,\eqno(13)
\end{equation*}%
where $B_{r}=\left\{ {\varphi \left( z\right) <\ln r}\right\} $.

We note, that $S^{\ast }$ - parabolic manifolds satisfy the condition (13).
In fact, if $\rho \left( z\right) $ is special exhaustion function, then $%
\left( {dd^{c}\rho }\right) ^{n}=0$ off a compact $K\subset \subset X$ \ so $%
\int_{B_{r}}{\left( {dd^{c}\rho }\right) ^{n}=}\int_{K}{\left( {dd^{c}\rho }%
\right) ^{n}=const\,,\,\,\,r\geqslant r_{0}}$. Hence, (13) holds.

If $X$ has a continuous plurisubharmonic exhaustion function satisfying the
condition (13), then every bounded above plurisubharmonic function on $X$ is
constant \cite{D}, so that this kind of manifolds are parabolic. In fact, a
more general result is also true.

\begin{theorem}
If on a Stein manifold $X$ of dimension $n$, there exist a plurisubharmonic
(not necessary continuous) exhaustion function $\varphi $ that satisfies, 
\begin{equation*}
\lim \inf_{r\rightarrow \infty }\frac{{\int_{B_{r}}{\left( {dd^{c}\varphi }%
\right) ^{n}}}}{{[\ln r]^{n}}}\,=\,0,\eqno(14).
\end{equation*}%
Then $X$ is parabolic.
\end{theorem}

\begin{proof}
Let's assume that $X$ satisfies the condition (14), but $X$ is not
parabolic. We take a sequence $1<r_{1}<r_{2}<...,r_{k}\rightarrow \infty ,$
such that 
\begin{equation*}
\lim_{r\rightarrow \infty }\frac{{\int_{B_{r_{k}}}{\left( {dd^{c}\varphi }%
\right) ^{n}}}}{{[\ln r_{k}]^{n}}}\,=\,0\eqno(15)
\end{equation*}%
Without loss a generality we can assume that the ball $B_{1}=\left\{ {%
\varphi \left( z\right) <0}\right\} \neq \emptyset $. Then according the
proposition 1 the P-measure $\omega ^{\ast }\left( {z,\overline{B}%
_{1},B_{r_{k}}}\right) $ decreases to $\omega ^{\ast }\left( {z,\overline{B}%
_{1}}\right) \neq -1$ as $k\rightarrow \infty $ . The function $\omega
^{\ast }\left( {z,\overline{B}_{1}}\right) $ is maximal, that is $%
(dd^{c}\omega ^{\ast })^{n}=0$ in $\,\,X\backslash B_{1}\,\,$ and is equal $%
\,\,-1\,\,$ on $\,\,\overline{B}_{1}$. Hence, by comparison principle of
Bedford-Taylor \cite{BT} we have: 
\begin{eqnarray*}
\int_{B_{r_{k}}}{\left[ {{\text{dd}}^{\text{c}}\omega ^{\ast }\left( {z,%
\overline{B}_{1},B_{r_{k}}}\right) }\right] ^{\text{n}}{\text{ }}} &\text{= }%
&{{}\int_{\overline{B}_{1}}{\left[ {{\text{dd}}^{\text{c}}\omega ^{\ast
}\left( {z,\overline{B}_{1},B_{r_{k}}}\right) }\right] ^{\text{n}}\geqslant }%
} \\
\int_{\overline{B}_{1}}{\left[ {{\text{dd}}^{\text{c}}\omega ^{\ast }\left( {%
z,\overline{B}_{1}}\right) }\right] ^{\text{n}}{\text{ }}} &\text{= }&{%
{}\alpha >0\,.{\text{ }}}
\end{eqnarray*}%
However, if we apply again the comparison principle to $\,\,\omega ^{\ast
}\left( {z,\overline{B}_{1},B_{r_{k}}}\right) \,\,$ and 
\begin{equation*}
\,\,w\left( z\right) =\frac{{\varphi \left( z\right) -\ln r_{k}}}{{\ln r_{k}}%
},
\end{equation*}%
then 
\begin{eqnarray*}
\frac{{\text{1}}}{{\left( {{\text{lnr}}_{\text{k}}}\right) ^{n}}}%
\int_{B_{r_{k}}}{\left[ {dd^{c}\varphi \left( z\right) }\right] ^{n}\,} &{=}%
&\int_{B_{r_{k}}}{\left[ {{\text{dd}}^{\text{c}}w\left( z\right) }\right] ^{%
\text{n}}\geqslant } \\
{\int_{B_{r_{k}}}{\left[ {{\text{dd}}^{\text{c}}\omega ^{\ast }\left( {z,%
\overline{B}_{1},B_{r_{k}}}\right) }\right] ^{\text{n}}}} &{\geqslant \alpha
>0\,.}&{{\text{ }}}
\end{eqnarray*}%
This contradiction proves the theorem.
\end{proof}

\bigskip

\textbf{3. Sibony-Wong manifolds. }

\bigskip

We next consider an important class of Stein manifolds (analytic sets) with
the Liouville property, which were introduced by Sibony - Wong \cite{SW}. To
describe these spaces we need to introduce some notation. For an $n$
dimensional closed subvariety $X$ of $\mathbb{C}_{w}^{N}$ let us denote, as
usual, by $\sigma ,$ the restriction of $\ln \left\vert w\right\vert $ on $%
X. $ Denoting the intersection of the $r$- ball in $\mathbb{C}^{N}$ with $X$
\ by $B_{r}=\left\{ {z\in X:\sigma \left( z\right) <\ln r}\right\} $ we can
describe Sibony - Wong class as those $X^{\prime }$s, such that 
\begin{equation*}
\lim_{r\rightarrow \infty }\frac{{vol\left( {B_{r}}\right) }}{{\ln r}}%
<\infty ,
\end{equation*}%
where the projective volume, $vol\left( {B_{r}}\right) \,$ is equal to $%
\frac{{H_{2n}\left( {B_{r}}\right) }}{{r^{2n}}}$ , $H_{2n}-$ the Hausdorff
measure (${\mathbb{R}}^{2n}$ -volume) of $B_{r}$. Sibony and Wong showed
that on such spaces any bounded holomorphic function is constant.

When $n=1,$ a special case of a result by Takegoshi \cite{T} states that if 
\begin{equation*}
\sup_{r}\frac{{vol\left( {B_{r}}\right) }}{{g}\left( {r}\right) }<\infty ,
\end{equation*}%
where ${g}:\mathbb{R}^{+}\rightarrow \mathbb{R}^{+}$ is a nondecreasing
continuous function such that 
\begin{equation*}
\int_{0}^{\infty }\frac{dr}{{g}\left( {r}\right) }=\infty ,
\end{equation*}%
then every negative smooth subharmonic function on\ $X$ reduces to a
constant , i.e. $X$ is parabolic.\emph{\ }

The proof of this proposition is based on the following estimation:%
\begin{equation*}
v\left( r\right) ^{2}\leq C{g}\left( {r}\right) \frac{d}{dr}\left( v\left(
r\right) \right) ,\,\,\forall \,v\in sh(X)\cap C^{1}(X),
\end{equation*}%
where $v\left( r\right) =\int_{B_{r}}dv\wedge d^{c}v\,\,$ and $\,\,C>0$ is a
constant. We note that if $v$ is an arbitrary subharmonic function we can
approximate it by smooth subharmonic functions $v_{j}\downarrow v,$ we
conclude that the above expression is also valid for arbitrary subharmonic
functions and hence the proof given in \cite{T} shows that such an $X$ is
parabolic. Taking $g(r)=\text{ln}r$, we see that 1-dimensional Sibony - Wong
manifolds are parabolic.

For $n>1$, taking into account that $vol\left( {B_{r}}\right) ={\int_{B_{r}}{%
\left( {dd^{c}\sigma }\right) ^{n}}}$ , by Wirtinger's theorem, we can
deduce from Theorem 3 above that $X$ is parabolic. Summarizing, we conclude
that Sibony-Wong manifolds are parabolic for any $n\in \mathbb{N}.$

In connection with Problem 2 of section 2 it will be of interest to
investigate $S^{\ast }$ -- parabolicity of Sibony-Wong manifolds. Affine
algebraic manifolds are among this class since their projective volume is
finite. Moreover they are $S^{\ast }$ -- parabolic as we have already seen.
On the other hand special exhaustion functions for $S^{\ast }$ -- parabolic
Sibony-Wong manifolds other than the algebraic ones; they can not be
asymptotically bigger than $\sigma \left( z\right) =\ln \left\vert
z\right\vert $ restricted to $X.$

\bigskip

\begin{theorem}
Let $X\subset \mathbb{C}^{N}$ be a closed submanifold and $\rho \left(
z\right)- $ a special exhaustion function on it. If 
\begin{equation*}
\underline{\lim }\frac{{\rho \left( z\right) }}{{\sigma \left( z\right) }}%
\geqslant \alpha >0,
\end{equation*}
then $X$ is\ an affine-algebraic set in $\mathbb{C}^{N}$.
\end{theorem}

\begin{proof}
Taking $C\rho $ instead $\rho $, if it is necessary, we can assume that,
there is some compact $K\subset \subset X$ such, that 
\begin{equation*}
\frac{{\rho \left( z\right) }}{{\sigma \left( z\right) }}\geqslant
1\,\,,\,\,z\in X\backslash K.
\end{equation*}%
Let $\sup_{K}\rho \left( z\right) =r_{0}$. Then $B_{r}=\left\{ {z\in X:\rho
\left( z\right) <\ln r}\right\} \,\,,\,\,r>r_{0},$ is not empty and open.
Hence, the closure $\overline{B}_{r}$ is not pluripolar. Therefore, the
extremal Green function 
\begin{equation*}
V_{\rho }\left( {z,\overline{B}_{r}}\right) =\sup \left\{ {u\left( z\right)
\in psh\left( X\right) :u|_{B_{r}}\leqslant 0,u\left( z\right) \leqslant
C_{u}+\rho \left( z\right) \,\,\forall z\in X}\right\}
\end{equation*}%
is locally bounded on $X$ (see \cite{Z3}). In the other hand, since $\rho
\left( z\right) \geqslant \sigma \left( z\right) $ outside of compact $K$,
then 
\begin{equation*}
V\left( {z,\overline{B}_{r}}\right) \leqslant V_{\rho }\left( {z,\overline{B}%
_{r}}\right) ,\text{where }V\left( {z,\overline{B}_{r}}\right) =V_{\sigma
}\left( {w,\overline{B}_{r}}\right) |_{X},
\end{equation*}%
\begin{equation*}
V_{\sigma }\left( {w,\overline{B}_{r}}\right) =\sup \left\{ {u\left(
w\right) \in psh\left( {\mathbb{C}^{N}}\right) :u|_{B_{r}}\leqslant 0,}\text{
}{u\left( w\right) \leqslant C_{u}+\ln \left\vert w\right\vert }\right\} .
\end{equation*}%
But the extremal function $V\left( {z,\overline{B}_{r}}\right) $ is locally
bounded on $X$ if and only if $X$ affine-algebraic \cite{SA3}. This
completes the proof.
\end{proof}

\bigskip

\textbf{4. Remark1:} Stoll in \cite{ST2} introduced and studied analytic
sets, for which the solution of the equation (in the notation of the above
section),%
\begin{equation*}
dd^{c}\omega _{R}\wedge \Psi =0\,\,,\,\,\omega _{R}|_{\partial
B_{0}}=-1\,,\,\,\,\omega _{R}|_{\partial B_{R}}=0\,,
\end{equation*}%
has the parabolic property, that $\omega _{R}\rightarrow -1\,\,,\text{for}%
\,\,R\nearrow \infty $, where $\Psi $ is close, positive $\left( {n-1\,,\,n-1%
}\right) $ form. Atakhanov \cite{ATA} called this kind of sets "parabolic
type" and proved that the sets which satisfy 
\begin{equation*}
\lim_{r\rightarrow \infty }\frac{{vol\left( {B_{r}}\right) }}{{\ln r}}=0
\end{equation*}%
are of this type. Moreover, he constructed \ Nevanlinna's equidistribution
theory for holomorphic maps $\,\,f:\mathrm{X}\rightarrow \mathrm{P}^{m}$. In
particular, on this kind of sets theorems of Picard, Nevanlinna, Valiron on
defect hyperplanes are true. \bigskip

\textbf{Remark 2}: In the literature there exits quite a number of
Liouville- type theorems for specific complex manifolds. However the
property that every bounded analytic function reduces to a constant need not
imply parabolicity, as is well known to people working in classification
theory of open Riemann surfaces. The simple example below illustrates this
point.

\bigskip

\textbf{Example:} Choose, on complex plane $\mathbb{C}_{z_{1}}$ a
subharmonic function $u$ with the property that $\left\{ {u\left( {z_{1}}%
\right) =-\infty }\right\} =\left\{ {0,1,\frac{1}{2},\frac{1}{3},...}%
\right\} .$ Let $\ w\left( {z_{1},z_{2}}\right) =u\left( {z_{1}}\right) +\ln
\left\vert {z_{2}}\right\vert .$ Then $w\in psh\left( {\mathbb{C}^{2}}%
\right) $, and the component $D$ of $\left\{ \left( {z_{1},z_{2}}\right) \in 
{\mathbb{C}}^{2}:w\left( {z_{1},z_{2}}\right) <0\right\} $ containing the
origin, being pseudoconvex, is a Stein manifold. Any bounded holomorphic
function on it is constant by the Liouville's theorem. However, the
plurisubharmonic function $w\left( {z_{1},z_{2}}\right) \neq const$ \ and is
bounded from above i.e. $D$ is not parabolic.

\bigskip

\begin{acknowledgement}
The authors wish to thank the referee for reading the manuscript carefully
and for valuable suggestions.
\end{acknowledgement}

\end{document}